\documentclass[12pt]{amsart}
\usepackage{amsmath,amssymb,amsbsy,amsfonts,amsthm,latexsym,amsopn,amstext,amsxtra,euscript,
amscd}
\usepackage{hyperref}
\usepackage{color}
\begin{document}

\title[Probabilistic interpretation of the M\"obius function ...]{Probabilistic interpretation \\
of the M\"obius function identity \\
and the Riemann Hypothesis}

\author[R. M.~Abrarov]{R. M.~Abrarov$^1$}
%\address{University of Toronto, Canada}
%\rabrarov@physics.utoronto.ca
\email{$^1$rabrarov@gmail.com}
\author[S. M.~Abrarov]{S. M.~Abrarov$^2$}
%\address{York University, Toronto, Canada}
%\email{abrarov@yorku.ca}
\email{$^2$absanj@gmail.com}

\date{February 8, 2010}

\begin{abstract}

We obtained the probabilities for the values of the M\"obius function for arbitrary numbers and found that the asymptotic densities of the squarefree integers among the odd and even numbers are $8/\pi^2$ and $4/\pi^2$, respectively. It is determined that statistics of successive outcomes of the M\"obius function for very large squarefree odd and even numbers behaves similar to statistics of heads and tails of two flipping coins. These preliminary results are giving arguments supporting the Riemann Hypothesis. Its plausibility is based on statistical phenomena for integers.

\end{abstract}

\maketitle

\section{\textbf{Probabilities in intervals}}

Consider in more details the identity (11) from Ref. \cite{RAbrarov1}
\begin{equation} \label{Eq_1} 
\mu \left(n\right)=-\sum _{i,j=1}^{\sqrt{n} }\mu _{i}  \mu _{j} \delta \left(\frac{n}{i\cdot j} \right),\qquad \qquad     n\ge 2,   
\end{equation}
where $\mu _{n} \equiv \mu \left(n\right)$  is the M\"obius function \cite{Edwards, Hardy, Weisstein1} and the delta function $\delta \left(\frac{n}{i\cdot j} \right)$ is 1 if $k$ divisible by $i\cdot j$ and 0 otherwise. Further, as a number we will imply a natural number.

In the identity \eqref{Eq_1} each term with the delta function $\delta \left(\frac{n}{i\cdot j} \right)$ for randomly chosen numbers $n$ can be interpreted as an event when the delta function became equal to 1. For the consecutive numbers such events occur through the period $i\cdot j$. It indicates that for the random numbers $n$ the frequency of event for the delta function $\delta \left(\frac{n}{i\cdot j} \right)$ is $\frac{1}{i\cdot j} $ and it can be considered as a probability of the event

\begin{equation} \label{Eq_2} 
\Pr \left(\delta \left(\frac{n}{i\cdot j} \right)\, =1\right)=\frac{1}{i\cdot j} .         
\end{equation} 

In \eqref{Eq_1} we have unchanging set of the terms only in the intervals between squares of numbers with consecutive nonzero values of the M\"obius function:
$$
\div 2^{2} \div 3^{2} \div 5^{2} \div 6^{2} \div 7^{2} \div 10^{2} \div 11^{2} \div 13^{2} \div 14^{2} \div 15^{2} \div\dots \, .
$$
Hence, the probability of the certain value of the M\"obius function is the dynamic function of numbers. For calculation of this probability we need to calculate from the formula \eqref{Eq_1} the probability of events leading to this specific value of the M\"obius function and to exclude the events leading to other values.

Let us illustrate above on simple examples, considering that randomly selected number falls into one of the intervals below.

i. \textbf{Interval $\left[2,2^{2} \right)$.}

$$
\mu \left(n\right)=-\delta \left(\frac{n}{1\cdot 1} \right) . 
$$
$$
\Pr \left(\mu _{n} =-1\right)=1\, ,\, \, \Pr \left(\mu _{n} =+1\right)=0\, ,\, \, \Pr \left(\mu _{n} =0\right)=0 .
$$

ii. \textbf{Interval $\left[2^{2} ,3^{2} \right)$.}

$$\mu \left(n\right)=-\delta \left(\frac{n}{1\cdot 1} \right)+2\delta \left(\frac{n}{1\cdot 2} \right)-\delta \left(\frac{n}{2\cdot 2} \right) .
$$
$\Pr \left(\mu _{n} =-1\right)=$Probability of odd numbers:
$$\Pr \left(\mu _{n} =-1\right)=1-\frac{1}{2} =\frac{1}{2} \left(1-2\frac{1}{2} +\frac{1}{2^{2} } \right)+\frac{1}{2} \left(1-\frac{1}{2^{2} } \right).
$$
$\Pr \left(\mu _{n} =+1\right)=$Probability of even numbers except numbers divisible by $2^{2} $:
$$
\Pr \left(\mu _{n} =+1\right)=\frac{1}{2} -\frac{1}{2^{2} } =-\frac{1}{2} \left(1-2\frac{1}{2} +\frac{1}{2^{2} } \right)+\frac{1}{2} \left(1-\frac{1}{2^{2} } \right) .
$$
$\Pr \left(\mu _{n} =0\right)=$Probability of numbers divisible by $2^{2} $.
$$
\Pr \left(\mu _{n} =0\right)=\frac{1}{2^{2} }  .
$$

iii. \textbf{Interval $\left[3^{2} ,\, \, 5^{2} \right)$.}

\begin{equation*}
\begin{aligned}
\mu \left(n\right)=&-\delta \left(\frac{n}{1\cdot 1} \right)+2\delta \left(\frac{n}{1\cdot 2} \right)-\delta \left(\frac{n}{2\cdot 2} \right) \\
&+2\delta \left(\frac{n}{1\cdot 3} \right)-2\delta \left(\frac{n}{2\cdot 3} \right)-\delta \left(\frac{n}{3\cdot 3} \right) .
\end{aligned}
\end{equation*}
$\Pr \left(\mu _{n} =-1\right)=$Probability of odd numbers except numbers divisible by 3:
\begin{equation*}
\begin{aligned}
\Pr \left(\mu _{n} =-1\right)&=1-\frac{1}{2} -\frac{1}{3} \left(1-\frac{1}{2} \right) \\
&=\frac{1}{2} \left(1-2\frac{1}{2} +\frac{1}{2^{2} } -2\frac{1}{3} +2\frac{1}{2\cdot 3} +\frac{1}{3^{2} } \right)\\
&+\frac{1}{2} \left(1-\frac{1}{2^{2} } -\frac{1}{3^{2} } \right) .
\end{aligned}
\end{equation*}
$\Pr \left(\mu _{n} =+1\right)=$Probability of even numbers and numbers divisible by 3 except numbers divisible by $2^{2} $ and $3^{2} $: 
\begin{equation*}
\begin{aligned}
\Pr \left(\mu _{n} =+1\right)=& \;\;\;\; \frac{1}{2} +\frac{1}{3} \left(1-\frac{1}{2} \right)-\left(\frac{1}{2^{2} } +\frac{1}{3^{2} } \right) \\
=&-\frac{1}{2} \left(1-2\frac{1}{2} +\frac{1}{2^{2} } -2\frac{1}{3} +2\frac{1}{2\cdot 3} +\frac{1}{3^{2} } \right) \\
&+\frac{1}{2} \left(1-\frac{1}{2^{2} } -\frac{1}{3^{2} } \right) .
\end{aligned}
\end{equation*}
$\Pr \left(\mu _{n} =0\right)=$Probability of numbers divisible by $2^{2} $ and $3^{2} $:
$$
\Pr \left(\mu _{n} =0\right)=\frac{1}{2^{2} } +\frac{1}{3^{2} } .
$$

Proceed to the general case continuing the same way. We have obvious relations: 

\begin{equation} \label{Eq_3} 
\Pr \left(\mu _{n} =+1\right)+\Pr \left(\mu _{n} =-1\right)+\Pr \left(\mu _{n} =0\right)=1 ,    
\end{equation} 
\begin{equation} \label{Eq_4} 
\Pr \left(\mu _{n} =0\right)=-\sum _{i=2}^{\sqrt{n} }\frac{\mu _{i} }{i^{2} }   ,        
\end{equation} 

\begin{equation} \label{Eq_5} 
\Pr \left(\left|\mu _{n} \right|=1\right)=\;\; \sum _{i=1}^{\sqrt{n} }\frac{\mu _{i} }{i^{2} }   .   
\end{equation}
It is not difficult to deduce that
\begin{equation} \label{Eq_6} 
\Pr \left(\mu _{n} =-1\right)=\;\; \frac{1}{2} \sum _{i,j=1}^{\sqrt{n} }\frac{\mu _{i} \mu _{j} }{i\cdot j}  +\frac{1}{2} \sum _{i=1}^{\sqrt{n} }\frac{\mu _{i} }{i^{2} }  =\;\;\frac{1}{2} m_{\sqrt{n} }^{2} +\frac{1}{2} \sum _{i=1}^{\sqrt{n} }\frac{\mu _{i} }{i^{2} }        ,    
\end{equation} 

\begin{equation} \label{Eq_7} 
\Pr \left(\mu _{n} =+1\right)=-\frac{1}{2} \sum _{i,j=1}^{\sqrt{n} }\frac{\mu _{i} \mu _{j} }{i\cdot j}  +\frac{1}{2} \sum _{i=1}^{\sqrt{n} }\frac{\mu _{i} }{i^{2} }  =-\frac{1}{2} m_{\sqrt{n} }^{2} +\frac{1}{2} \sum _{i=1}^{\sqrt{n} }\frac{\mu _{i} }{i^{2} }   ,    
\end{equation} 

\begin{equation} \label{Eq_8} 
\Delta \Pr \left(n\right)=\Pr \left(\mu _{n} =-1\right)-\Pr \left(\mu _{n} =+1\right)=m_{\sqrt{n} }^{2}  ,      
\end{equation} 
where $m_{n} =\sum _{i=1}^{n}\frac{\mu _{i} }{i}  $.

From \eqref{Eq_8} we see that the difference between the probabilities in the above intervals as a function of $n$ for the negative and positive values of the M\"obius function tend to zero as $m_{\sqrt{n} }^{2} $. Such advancement provides an opportunity to apply the Denjoy's proposal \cite{Denjoy, Edwards} and induction procedure \cite{RAbrarov2} for verification of the Riemann Hypothesis \cite{Bombieri, Conrey, Edwards, Wikipedia}.

The Denjoy's proposal can be briefly described as follows. Suppose for squarefree numbers (integer $k$ is squarefree if $\mu _{k} \ne 0$ \cite{Weisstein2}) that $\{ \mu _{k} \} $ is a sequence of random and independent events with symmetrical distribution
\begin{equation} \label{Eq_9} 
\Pr \left(\mu _{k} =+1\right)=\Pr \left(\mu _{k} =-1\right)=\frac{1}{2} .       
\end{equation} 

Let us start summation $S_{n} =\sum _{k=m}^{n}\mu _{k}  $ from some integer $m$ (it can be, for example, a very large number). From the de Moivre-Laplace limit theorem follows that
\begin{equation} \label{Eq_10} 
\mathop{\lim }\limits_{n\to \infty } \Pr \left(\left|S_{n} \right|\le c\sqrt{n} \right)=\frac{1}{\sqrt{2\pi } } \int _{-c}^{c}\exp \left(-\frac{x^{2} }{2} \right)dx  .     
\end{equation} 
Since the right hand side of equation above tends to 1 as $c\to \infty $, for every  $\varepsilon>$0 we have
\begin{equation} \label{Eq_11} 
\mathop{\lim }\limits_{n\to \infty } \Pr \left(\left|S_{n} \right|<n^{\frac{1}{2} +\varepsilon } \right)=1.        
\end{equation} 
If values $\{ \mu _{k} \} $ of the M\"obius function for squarefree numbers would satisfy the above conditions, then the Riemann hypothesis is true with probability one \cite{Edwards}.

\section{\textbf{Probabilities for arbitrary numbers}}

However, it is still unclear whether the Denjoy's proposal is applicable here since it is unknown how independent the successive evaluations of the M\"obius function.

Let us try to find from the probabilities in intervals \eqref{Eq_3}-\eqref{Eq_8} the probabilities for arbitrary numbers. For odd numbers we proceed from the basic identity \eqref{Eq_1} excluding all delta functions with even divisors
\begin{equation} \label{Eq_12} 
\begin{aligned}
\mu \left(2n-1\right)&=-\sum _{i,j=1}^{\sqrt{2n-1} }\mu _{i}  \mu _{j} \delta \left(\frac{2n-1}{i\cdot j} \right) \\
&=-\sum _{i,j=1\left(odd\right)}^{\sqrt{2n-1} }\mu _{i}  \mu _{j} \delta \left(\frac{2n-1}{i\cdot j} \right) ,   
\end{aligned}
\end{equation} 
where (odd) denotes the summation over odd indices only. Consequently, by analogy with \eqref{Eq_3}-\eqref{Eq_8} we have
\begin{equation} \label{Eq_13} 
\begin{aligned}
\Pr \left(\mu _{2n-1} =-1\right)&=\;\; \frac{1}{2} \sum _{i,j=1\left(odd\right)\, }^{\sqrt{2n-1} }\frac{\mu _{i} \mu _{j} }{i\cdot j}  +\frac{1}{2} \sum _{i=1\left(odd\right)}^{\sqrt{2n-1} }\frac{\mu _{i} }{i^{2} } \\
&  =\;\; \frac{1}{2} \left(m_{\sqrt{2n-1} }^{odd} \right)^{2} +\frac{1}{2} \sum _{i=1\left(odd\right)}^{\sqrt{2n-1} }\frac{\mu _{i} }{i^{2} }  , 
\end{aligned}
\end{equation} 

\begin{equation} \label{Eq_14}
\begin{aligned}
\Pr \left(\mu _{2n-1} =+1\right)&=-\frac{1}{2} \sum _{i,j=1\left(odd\right)\, }^{\sqrt{2n-1} }\frac{\mu _{i} \mu _{j} }{i\cdot j}  +\frac{1}{2} \sum _{i=1\left(odd\right)}^{\sqrt{2n-1} }\frac{\mu _{i} }{i^{2} } \\
& =-\frac{1}{2} \left(m_{\sqrt{2n-1} }^{odd} \right)^{2} +\frac{1}{2} \sum _{i=1\left(odd\right)}^{\sqrt{2n-1} }\frac{\mu _{i} }{i^{2} }   ,
\end{aligned}
\end{equation} 

\begin{equation} \label{Eq_15} 
\Pr \left(\mu _{2n-1} =0\right)=-\sum _{i=3\left(odd\right)}^{\sqrt{2n-1} }\frac{\mu _{i} }{i^{2} }    ,        
\end{equation} 

\begin{equation} \label{Eq_16} 
\Pr \left(\left|\mu _{2n-1} \right|=1\right)=\;\; \sum _{i=1\left(odd\right)}^{\sqrt{2n-1} }\frac{\mu _{i} }{i^{2} }    ,        
\end{equation} 
where $m_{n}^{odd} \equiv \sum _{i=1\left(odd\right)\, }^{n}\frac{\mu _{i} }{i}  $. 

Probability of squarefree numbers among odd numbers tends to the limit (asymptotic density)
\begin{equation} \label{Eq_17} 
\mathop{\lim }\limits_{n\to \infty } \Pr \left(\left|\mu _{2n-1} \right|=1\right)=\sum _{i=1\left(odd\right)}^{\infty }\frac{\mu _{i} }{i^{2} } = \frac{8}{\pi ^{2} }  .       
\end{equation} 
It is interesting to note that we can go further and exclude from identity \eqref{Eq_12} the delta functions with divisors 3, then 5,  ... etc. ultimately getting the trivial result
$$
\Pr \left(\mu _{prime} =-1\right)=1.
$$
For even numbers we can use the obtained equations \eqref{Eq_6}-\eqref{Eq_7}, which are, as we see, the arithmetic average for the probabilities of odd and even numbers 

\begin{equation} \label{Eq_18} 
\Pr \left(\mu _{n}^{even} =\pm 1,\, \, 0\right)=2\Pr \left(\mu _{n} =\pm 1,\, \, 0\right)-\Pr \left(\mu _{n}^{odd} =\pm 1,\, \, 0\right)  ,    
\end{equation} 

\begin{equation} \label{Eq_19} 
\begin{aligned}
\Pr \left(\mu _{2n} =-1\right)&=\;\; \frac{1}{2} \left(2\sum _{i,j=1\, }^{\sqrt{2n} }\frac{\mu _{i} \mu _{j} }{i\cdot j}  -\sum _{i,j=1\left(odd\right)\, }^{\sqrt{2n} }\frac{\mu _{i} \mu _{j} }{i\cdot j}  \right) \\
&\;\;\; +\frac{1}{2} \left(2\sum _{i=1}^{\sqrt{2n} }\frac{\mu _{i} }{i^{2} }  -\sum _{i=1\left(odd\right)}^{\sqrt{2n} }\frac{\mu _{i} }{i^{2} }  \right)  ,  
\end{aligned}
\end{equation} 

\begin{equation} \label{Eq_20} 
\begin{aligned}
\Pr \left(\mu _{2n} =+1\right)=&-\frac{1}{2} \left(2\sum _{i,j=1\, }^{\sqrt{2n} }\frac{\mu _{i} \mu _{j} }{i\cdot j}  -\sum _{i,j=1\left(odd\right)\, }^{\sqrt{2n} }\frac{\mu _{i} \mu _{j} }{i\cdot j}  \right) \\
&+\frac{1}{2} \left(2\sum _{i=1}^{\sqrt{2n} }\frac{\mu _{i} }{i^{2} }  -\sum _{i=1\left(odd\right)}^{\sqrt{2n} }\frac{\mu _{i} }{i^{2} }  \right)  ,  
\end{aligned}
\end{equation} 

\begin{equation} \label{Eq_21} 
\Pr \left(\mu _{2n} =0\right)=-\left(2\sum _{i=2}^{\sqrt{2n} }\frac{\mu _{i} }{i^{2} }  -\sum _{i=3\left(odd\right)}^{\sqrt{2n} }\frac{\mu _{i} }{i^{2} }  \right)   ,      
\end{equation} 

\begin{equation} \label{Eq_22} 
\Pr \left(\left|\mu _{2n} \right|=1\right)=2\sum _{i=1}^{\sqrt{2n} }\frac{\mu _{i} }{i^{2} }  -\sum _{i=1\left(odd\right)}^{\sqrt{2n} }\frac{\mu _{i} }{i^{2} }     ,      
\end{equation} 

\begin{equation} \label{Eq_23} 
\mathop{\lim }\limits_{n\to \infty } \Pr \left(\left|\mu _{2n} \right|=1\right)=2\sum _{i=2}^{\infty }\frac{\mu _{i} }{i^{2} }  -\sum _{i=1\left(odd\right)}^{\infty }\frac{\mu _{i} }{i^{2} }  =\frac{4}{\pi ^{2} }  .      
\end{equation} 
Similar to \eqref{Eq_8}, we see that difference between probabilities for negative and positive values of the M\"obius function for odd and even numbers tend to zero when $n\to \infty $:

\begin{equation} \label{Eq_24}
\begin{aligned}
\Delta \Pr \left(2n-1\right)&=\Pr \left(\mu _{2n-1} =-1\right)-\Pr \left(\mu _{2n-1} =+1\right) \\
&=\left(m_{\sqrt{2n-1} }^{odd} \right)^{2},    
\end{aligned}
\end{equation} 

\begin{equation} \label{Eq_25} 
\begin{aligned}
\Delta \Pr \left(2n\right)&=\Pr \left(\mu _{2n} =-1\right)-\Pr \left(\mu _{2n} =+1\right) \\
&=2m_{\sqrt{2n} }^{2} -\left(m_{\sqrt{2n} }^{odd} \right)^{2}  .   
\end{aligned}
\end{equation} 
Now we have the probabilities for arbitrary numbers and can to come to some conclusions.

\section{\textbf{Conclusions}}

All probabilities for the values of the M\"obius function show oscillatory behaviour with tendency to approach the certain limits. We obtained the expected result that among odd numbers the squarefree integers occur more frequently than for even numbers and calculated their asymptotic densities \eqref{Eq_17} and \eqref{Eq_23}.  For very large integers, where oscillation of the probabilities is much less than absolute values of their limits, the squarefree numbers occur with almost equal probabilities (symmetrically) for positive and negative values of the M\"obius function. Moreover, they occur symmetrically regardless whether it is odd or even number and what number was before; the outcome of the earlier event for extremely large squarefree number does not affect the next event. Therefore these events are independent. Intuitively it is well explainable since an appearance of very large squarefree number involves an enormous quantity of the delta functions generating the M\"obius function outcome \eqref{Eq_1}. Hence, we have intermixing of colossal quantity of the events \eqref{Eq_2} leading to intense randomization and chaos. In this sense it is related to the process of the Brownian motion - the delta functions in \eqref{Eq_1} are acting as a chaotic process of collisions of the fluid molecules causing the chaotic movement of the Brownian particle. There is another association as well - if we exclude statistics for squareful numbers (numbers with $\mu_n=0$ \cite{Weisstein2}), we can say that very large odd and even numbers behave as two unbiased tossing coins with only difference that instead of heads and tails they produce positive and negative values of the M\"obius function. Consequently, all conditions of Denjoy's proposal \cite{Denjoy} can be applied. Some dominance of the probability for negative values of the M\"obius function vanishing at $n\to \infty $ \eqref{Eq_8} leads to the systematic shift of the Mertens function \cite{Edwards, Hardy, Wikipedia} of the order $n\cdot m_{\sqrt{n} }^{2} $ towards the negative area. Such tendency can be observed even in the early evolution of the Mertens function from the Mathematica demonstration \cite{WolframDemo}. It is possible to show by induction \cite{RAbrarov2} that this shift remains at the same order as oscillations of the Mertens function due to random character of occurrences of squarefree numbers.

These preliminary results are giving arguments supporting the Riemann Hypothesis. Its plausibility is based on statistical phenomena for integers.
\\

\end{document}